\documentstyle[11pt,amssymb,amsfonts]{article}

\begin{document}
\newcommand{\p}{\parallel }
\makeatletter \makeatother
\newtheorem{th}{Theorem}[section]
\newtheorem{lem}{Lemma}[section]
\newtheorem{de}{Definition}[section]
\newtheorem{rem}{Remark}[section]
\newtheorem{cor}{Corollary}[section]
\renewcommand{\theequation}{\thesection.\arabic {equation}}

\title{{\bf The noncommutative infinitesimal equivariant index formula: part II}}

\author{ Yong Wang }

\date{}
\maketitle

\begin{abstract}
In this paper, we prove that infinitesimal equivariant Chern-Connes
characters are well-defined. We decompose an equivariant index as a pairing of infinitesimal equivariant Chern-Connes
characters with the Chern character of an idempotent matrix. We compute the limit of infinitesimal equivariant Chern-Connes
characters when the time goes to zero by using the Getzler symbol calculus and then extend these theorems to the family case. We also prove that infinitesimal equivariant eta cochains are well-defined and prove the noncommutative infinitesimal equivariant index formula for manifolds with boundary.\\

\\


\end{abstract}

\section{Introduction}

    \quad The Atiyah-Bott-Segal-Singer index formula is a
generalization of the Atiyah-Singer index theorem to manifolds admitting group actions.
In [BV1], [LYZ], [PW], various heat kernel proofs of the equivariant index theorem have been given and each method
has its own advantage. For manifolds with boundary, the equivariant extension of the Atiyah-Patodi-Singer index theorem was given by
Donnelly in [Do]. In the equivariant Atiyah-Patodi-Singer index theorem, the equivariant eta invariant appears and the regularity of the equivariant eta invariant was proved by Zhang in [Zh].
An infinitesimal version of the equivariant index formula was established in [BV2] and a direct heat kernel proof was given by Bismut in [Bi].
The infinitesimal equivariant index formula for manifolds with boundary
was established in [Go] with the introduction of the infinitesimal
equivariant eta invariant.\\
\indent The counterpart of the index formula in the noncommutative geometry is the computation of the Chern-Connes character [Co], [JLO], [GS].
The JLO character was computed in [CM] and [BF] by using the Getzler symbol calculus in [Ge2].
In [Az], [CH] and [PW], these authors gave the computations of the equivariant JLO characters
associated to a G-equivariant
$\theta$-summable Fredholm module.
In [Wa], we defined the truncated infinitesimal equivariant Chern-Connes
characters and computed the limit of the truncated infinitesimal equivariant Chern-Connes
characters when the time goes to zero.\\
\indent Compared with [Wa], there are several improvements in the present paper. In (2.2) in [Wa], we defined truncated infinitesimal equivariant Chern-Connes
characters. It is only well-defined when it is a polynomial of Lie algebra elements. In this paper, we drop off
the truncated order $J$ (see (2.2)) and this consequently requires much better estimates (see Lemma 2.2). As in [GS], we decompose an equivariant index as a pairing of infinitesimal equivariant Chern-Connes characters with the Chern character of an idempotent matrix. Compared with Corollary
2.13 in [Wa], we drop off the limit on the right hand side of Corollary 2.13. Next we compute the limit of infinitesimal equivariant Chern-Connes
characters when the time goes to zero by using the Getzler symbol calculus. Since we have dropped off the truncated order, (2.15) in [Wa] does not hold for our infinitesimal equivariant Chern-Connes characters. So we can not directly apply the method of
Theorem 2.12 in [Wa]. Instead, we first apply the Getzler symbol calculus to prove the existence of the limit of infinitesimal equivariant Chern-Connes characters when time goes to zero (Theorem 3.9) and then use Theorem 2.12 in [Wa] to get the result. On the direction, in Section 3 in [Wa], we define the truncated infinitesimal equivariant eta cochains. Again in this paper we drop off the truncated order and then give a proof of the regularity at zero of infinitesimal equivariant eta cochains by using the method in [PW]. That is, we prove that (3.5) in [Wa] holds for any $k$. This allows us to establish the noncommutative infinitesimal equivariant index formula for manifolds with boundary (see Theorem 4.9). In this paper, we also define family infinitesimal equivariant Chern-Connes characters and give the family generalization of the above theorems which does not appear in [Wa].\\
\indent This paper is organized as follows: In Section 2, we prove that infinitesimal equivariant Chern-Connes
characters are well-defined. Then we decompose the equivariant index as a pairing of infinitesimal equivariant Chern-Connes
characters with the Chern character of an idempotent matrix. In Section 3, We compute the limit of infinitesimal equivariant Chern-Connes
characters when the time goes to zero by using the Getzler symbol calculus. In Section 4, we prove that infinitesimal equivariant eta cochains are well-defined
and prove the noncommutative infinitesimal equivariant index formula for manifolds with boundary.
In Section 5, we extend results in Section 2, 3 to the family case.\\

\section{ The infinitesimal equivariant JLO cocycle and the index pairing}

\quad Let $M$ be a compact
oriented even dimensional Riemannian manifold without boundary with a
fixed spin structure and $S$ be the bundle of spinors on $M$. Denote
by $D$ the associated Dirac operator on $H=L^2(M;S)$, the Hilbert
space of $L^2$-sections of the bundle $S$. Let $c(df): S\rightarrow
S$ denote the Clifford action with $f\in C^{\infty}(M)$. Suppose
that $G$ is a compact connected Lie group acting on $M$ by
orientation-preserving isometries preserving the spin structure and $\mathfrak{g}$ is the Lie algebra of $G$. Then $G$ commutes with the Dirac operator. For $X\in\mathfrak{g}$, let $X_M(p)=\frac{d}{dt}|_{t=0}e^{-tX}p$ be the Killing
field induced by $X$. Let $c(X)$ denote the Clifford action by $X_M$, and ${\mathfrak{L}}_X$ denote the Lie derivative respectively. Define $\mathfrak{g}$-equivariant modifications of $D$ and $D^2$ for $X\in {\mathfrak{g}}$ as follows:\\
$$D_X:=D-\frac{1}{4}c(X);~~H_X:=D^2_{-X}+{\mathfrak{L}}_X=(D+\frac{1}{4}c(X))^2+{\mathfrak{L}}_X.\eqno(2.1)$$
Then $H_X$ is the equivariant Bismut Laplacian. Let ${\mathbb{C}}[{\mathfrak{g}}^*]$ denote the space of formal power series in $X\in {\mathfrak{g}}$ and $\psi_t$ be the rescaling operator on ${\mathbb{C}}[{\mathfrak{g}}^*]$ which is defined by $X\rightarrow \frac{X}{t}$ for $t>0$.\\
 \indent Let
$$A=C_G^{\infty}(M)=\{f\in C^{\infty}(M)|f(g\cdot x)=f(x), g\in~ G, x\in~ M\},$$
 then the data $(A,H,D+\frac{1}{4}c(X),G)$ defines a non selfadjoint perturbation of finitely summable
(hence $\theta$-summable) equivariant unbounded Fredholm module
$(A,H,D,G)$ in the sense of [KL] (for details, see [CH] and [KL]).
For $(A,H,D+\frac{1}{4}c(X),G)$,
  {\bf the infinitesimal equivariant JLO cochain} ${{\bf
{\rm ch}}}^{
2k}(D,X)$ can be defined by the formula:\\

$${{\bf {\rm ch}}^{2k}}(D,X)(f^0,\cdots ,
f^{2k}):
=\int_{\triangle_{2k}}{\rm
Str}\left[ e^{-{
\mathfrak{L}}_X}f^0e^{-\sigma_0(D+\frac{1}{4}c(X))^2}c(df^1)\right.$$
$$\left.\cdot e^{-\sigma_1(D+\frac{1}{4}c(X))^2}\cdots c(df^{2k})e^{-\sigma_{2k}(D+\frac{1}{4}c(X))^2}\right]d{\rm Vol}_{\Delta_{2k}},\eqno(2.2)$$
 where
$\triangle_{2k}=\{(\sigma_0,\cdots,\sigma_{2k})|~\sigma_0+\cdots \sigma_{2k}=1\}$ is the $2k$-simplex.
For an integer $J\geq 0$, denote by ${\mathbb{C}}[{\mathfrak{g}}^*]_J$ the space of polynomials of degree
$\leq J$ in $X\in {\mathfrak{g}}$ and let $(\cdot)_J:~ {\mathbb{C}}[{\mathfrak{g}}^*]\rightarrow {\mathbb{C}}[{\mathfrak{g}}^*]_J$ be the natural projection. Fix basis $e_1,\cdots,e_n$ of ${\mathfrak{g}}$ and let $X=x_1e_1+\cdots x_ne_n.$ A $J$-degree polynomial on $X$ is namely a $J$-degree polynomial on $x_1,\cdots,x_n.$ Now we prove that
${{\bf {\rm ch}}^{2k}}(D,X)(f^0,\cdots ,
f^{2k})$ is well-defined.\\
\indent Let $H$ be a Hilbert space. For
$q\geq0$, denote by $||.||_q$ the Schatten $p$-norm on the Schatten
ideal
$L^p$. Let $L(H)$ denote the Banach algebra of bounded operators on $H$.\\

\noindent {\bf Lemma 2.1}~([Si]){\it~~(i)~~${\rm Tr}(AB)={\rm
Tr}(BA)$, for $A,~B\in L(H)$ and $AB, ~BA\in
L^1$.\\
~~(ii)~~For $A\in L^1,$ we have $|{\rm Tr}(A)|\leq ||A||_1$,
$||A||\leq ||A||_1$.\\
~~(iii)~~For $A\in L^q$ and $B\in L(H)$, we have: $||AB||_q\leq
||B||||A||_q$, $||BA||_q\leq ||B||||A||_q$.\\
 ~~(iv)~(H\"{o}lder
Inequality)~~If $\frac{1}{r}=\frac{1}{p}+\frac{1}{q},~p,q,r>0,~A\in
L^p,~B\in
L^q,$ then $AB\in L^r$ and $||AB||_r\leq ||A||_p||B||_q$.}\\

\indent Let $H_X=D^2+F_X$, where $F_X$ is a first order differential operator with degree $\geq 1$ coefficients depending on $X$.\\

\noindent{\bf Lemma 2.2} {\it For any $1\geq u> 0$,~$t>0$, we have:}\\
$$||e^{-utH_X}||_{u^{-1}}\leq 2e^{\frac{t}{2}}\{1+[||(1+D^2)^{-\frac{1}{2}}F_X||^{2}e^{-1}\pi ut]^{\frac{1}{2}}\}e^{||(1+D^2)^{-\frac{1}{2}}F_X||^{2}e^{-1}\pi ut}
 ({\rm
tr}[e^{-\frac{tD^2}{2}}])^u.
\eqno(2.3)$$

\noindent {\bf Proof.} By the Duhamel principle, it is that
$$||e^{-utH_X}||_{u^{-1}}=||
\sum_{m\geq
0}(-ut)^m\int_{\triangle_m}e^{-v_0utD^2}F_Xe^{-v_1utD^2}$$
$$\cdot F_X\cdots e^{-v_{m-1}utD^2}F_Xe^{-v_mutD^2}dv||_{u^{-1}}.\eqno(2.4)$$
Also $||
(-ut)^m\int_{\triangle_m}e^{-v_0utD^2}F_Xe^{-v_1utD^2} F_X\cdots e^{-v_{m-1}utD^2}F_Xe^{-v_mutD^2}dv||_{u^{-1}}$ is continuous and bounded
by (2.7) in [Wa]. By the measure of the boundary of ${\triangle_m}$ being zero,
we can estimate (2.4) in the interior of ${\triangle_m}$, that is $v_j>0$. It holds that
$$||e^{-\frac{v_j}{2}utD^2}F_X||\leq (v_jut)^{-\frac{1}{2}}e^{-\frac{1-v_jut}{2}}||(1+D^2)^{-\frac{1}{2}}F_X||,\eqno(2.5)$$
where we use that $F_{X}$ is a first order differential
operator and the equality
$${\rm
sup}\{(1+x)^{\frac{l}{2}}e^{-\frac{utx}{2}}\}=(ut)^{-\frac{l}{2}}e^{-\frac{l-ut}{2}}.\eqno(2.6)$$
By the H\"{o}lder inequality, (2.4) and (2.5), the conditions that $0<u\leq 1$ and $v_0+\cdots+v_{m-1}\leq 1$ , we have
$$||e^{-utH_X}||_{u^{-1}}\leq e^{\frac{t}{2}}\sum_{m\geq 0} e^{-\frac{m}{2}}(ut)^{\frac{m}{2}}||(1+D^2)^{-\frac{1}{2}}F_X||^m
\int_{\triangle_m}v_0^{-\frac{1}{2}}\cdots v_{m-1}^{-\frac{1}{2}}dv({\rm
tr}[e^{-\frac{tD^2}{2}}])^u.\eqno(2.7)$$
It holds that (see line 7 in [BC, P.21])
$$\int_{\triangle_m}v_0^{-\frac{1}{2}}\cdots v_{m-1}^{-\frac{1}{2}}dv=\frac{\pi^{\frac{m}{2}}}{\frac{m}{2}\Gamma(\frac{m+1}{2})}.\eqno(2.8)$$
By $\Gamma(x+1)=x\Gamma(x),$ $\Gamma(n)=(n-1)!$ and $\Gamma(\frac{1}{2})=\sqrt{\pi},$ then
$\Gamma(\frac{m+1}{2})=(\frac{m-1}{2})!$ when $m$ is odd; $\Gamma(\frac{m+1}{2})=\frac{(m-1)!!\sqrt{\pi}}{2^{\frac{m}{2}}}$ when $m$ is even.
By (2.8) and
$${\rm lim}_{m\rightarrow +\infty}\frac{(2m-1)!!}{(2m)!!}=0,\eqno(2.9)$$
we know that the series (2.7) is absolutely convergent.
When $m$ is odd, then
$$\frac{\pi^{\frac{m}{2}}}{\frac{m}{2}\Gamma(\frac{m+1}{2})}\leq \frac{2\pi^{\frac{m}{2}}}{(\frac{m+1}{2})!}.\eqno(2.10)$$
When $m$ is even, then
$$\frac{\pi^{\frac{m}{2}}}{\frac{m}{2}\Gamma(\frac{m+1}{2})}\leq \frac{2\pi^{\frac{m}{2}}}{(\frac{m}{2})!}.\eqno(2.11)$$
By (2.7), (2.8), (2.10) and (2.11), we have
$$||e^{-utH_X}||_{u^{-1}}\leq 2e^{\frac{t}{2}}\left[\sum_{m~ {\rm even}}
\frac{(||(1+D^2)^{-\frac{1}{2}}F_X||^2e^{-1}\pi ut)^{\frac{m}{2}}}{(\frac{m}{2})!}\right.$$
$$\left.
+\sum_{m~ {\rm odd}}
\frac{(||(1+D^2)^{-\frac{1}{2}}F_X||^2e^{-1}\pi ut)^{\frac{m}{2}}}{(\frac{m+1}{2})!}\right]
({\rm
tr}[e^{-\frac{tD^2}{2}}])^u.\eqno(2.12)$$
Therefore, (2.3) can be obtained.~~$\Box$\\

\indent By (2.2), (2.3) and the H\"{o}lder inequality as well as ${\rm Vol}_{\triangle_{2k}}=\frac{1}{(2k)!}$, for $t=1$ and $\sigma_l\leq 1$,
we get
$$|{{\bf {\rm ch}}^{2k}}(D,X)(f^0,\cdots ,
f^{2k})|\leq \frac{1}{(2k)!}||f^0||(\prod_{j=1}^{2k}||df^j||)$$
$$\cdot
[2e^{\frac{1}{2}}
(1+(||(1+D^2)^{-\frac{1}{2}}F_X||^2e^{-1}\pi)^{\frac{1}{2}})]^{2k+1}
e^{||(1+D^2)^{-\frac{1}{2}}F_X||^{2}e^{-1}\pi}
({\rm
tr}[e^{-\frac{D^2}{2}}]).\eqno(2.13)$$
Thus, ${{\bf {\rm ch}}^{2k}}(D,X)$ is well-defined. Recall that an even cochain $\{\Phi_{2n}\}$ is called
entire if $\sum_n||\Phi_{2n}||n!z^n$ is entire, where $||\Phi||:={\rm sup}_{||f^j||_1\leq 1}\{|\Phi(f^0,f^1,\cdots,f^{2k})|\}$.
By (2.13), then $\{{{\bf {\rm ch}}^{2k}}(D,X)\}$ is an entire cochain. Let $p\in M_r({\mathbb{C}}^{\infty}(M))$ and $p=p^2=p^*$ and
$p(gx)=p(x)$. Define
the Chern character of $p$ by (see [GS])
$${\rm ch}(p):={\rm Tr}(p)+\sum_{l}(-1)^l\frac{(2l)!}{2\cdot l!}{\rm Tr}(2p-1,p,\cdots,p)_{2l}.\eqno(2.14)$$
By (2.13), $\left<{\rm ch}^*(D,X),{\rm ch}(p)\right>$ is convergent.
Similarly to Theorem A in [GS], we have\\

\noindent {\bf Proposition 2.3} (1) {\it  The infinitesimal equivariant Chern-Connes character is closed:}
$$(B+b)({\rm ch}^*(D,X))=0.\eqno(2.15)$$
(2)~ {\it Let $D_{\tau}=D+\tau V$ and $D_{-X,\tau}=D_{-X}+\tau V$ and $V$ is a bounded operator which commutes with $e^{-X}$, then there exists a cochain ${\rm ch}^*(D_{\tau},X,V)$ such that}
$$\frac{d}{d\tau}{\rm ch}^*(D_{\tau},X)=-(B+b){\rm ch}^*(D_{\tau},X,V).\eqno(2.16)$$\\

\indent By the Serre-Swan theorem, we denote the vector bundle over $M$ with the fibre $p(x)({\mathbb{C}}^r)$ at $x\in M$ by ${{\rm Im}p}$.
 Let $D_{{\rm Im}p}$ be the Dirac operator twisted by the bundle ${{\rm Im}p}$.
By Proposition 2.3, $(B+b){\rm ch}(p)=0$ and Proposition 8.11 in [BGV], we have by taking $V=(2p-1)[D,p]$ that (see Section 3 in [GS])\\

\noindent{\bf Theorem 2.4} {\it The following index formula holds}
$${\rm Ind}_{e^{-X}}(D_{{\rm Im}p,+}) =\left<{\rm ch}^*(D,X),{\rm ch}(p)\right>.\eqno(2.17)$$\\

In Theorem 2.4, $X$ is unnecessarily small.\\

\section{The computations of infinitesimal equivariant Chern-Connes characters}

\quad In this section, we will compute infinitesimal equivariant Chern-Connes characters by Theorem 2.12 in [Wa] and the Getzler symbol
calculus in [Ge2] and [BF]. Recall the Getzler symbol calculus in [Ge2] and [BF]. Let $E$ be a vector bundle over the compact manifold $M$ and $\pi: T^*M\rightarrow M$ be the natural map and $E^0=\pi^*({\rm Hom}(E,E))$ be the pull-back of the bundle $\pi^*({\rm Hom}(E,E))$ to a bundle over $T^*M$.\\

\noindent{\bf Definition 3.1} A section $p\in E^0$ is called a symbol of order $l$ if for every multi-index $\alpha$ and $\beta$ we have the estimates:
$$||\partial^\alpha_x\partial^\beta_\xi p(x,\xi)||\leq C_{\alpha\beta}(1+|\xi|)^{m-|\beta|}.\eqno(3.1)$$
We denote by $\Sigma^l(E)$ the symbols of order $l$.\\

By the representative theorem of the Clifford algebra $Cl(T^*M) \simeq {\rm Hom}(S(TM))$ and the isomorphism $Cl(T^*M)\simeq \wedge(T^*M)$, note
a map $\overline{\sigma}$ defined by $$\overline{\sigma}:~{\rm Hom}(S(TM)\otimes E)\simeq {\rm Hom }E\otimes Cl(T^*M)\simeq {\rm Hom }E\otimes \wedge(T^*M),\eqno(3.2)$$
and $\overline{\theta}$ is the inverse of $\overline{\sigma}$. Let ${\mathbb{L}}=\pi^*({\rm Hom(E)}\otimes \wedge(T^*M))
\otimes {\mathbb{C}}[\mathfrak{g}^*]$ and $X\in {\mathfrak{g}}$. \\

\noindent{\bf Definition 3.2} A section $p\in {\mathbb{L}}$ is called a $s$-symbol of order $l$ if
$$p=\sum^{{\rm dim}M}_{j=0}(\sum_{|\alpha|\geq 0}p_{j,\alpha}X^\alpha)\otimes \omega_j,\eqno(3.3)$$
where~$\omega_j\in \Omega^j(M)$,~$p_{j,\alpha}\in
\Sigma^{l-j-2|\alpha|}(E) $ and $\sum_{|\alpha|\geq 0}||p_{j,\alpha}(x,\xi)|||X^\alpha|$ is convergent.
We denote the collection of $s$-symbol of order $l$ by $S\Sigma^l(E,X)$.\\

Let $x_0$ be a fixed point in $M$ and $T_{x_0}M$ be the tangent space and ${\rm exp}$ be the exponential map respectively. Let
$h$ be a function that is identically one in a neighborhood of the diagonal of $M\times M$ such that the exponential map is a diffeomorphism
on the support of $h$. Let $(x_0,x)\in {\rm supp}(h)$. Let $\tau(x_0,x):(S(TM)\otimes E)_{x_0}\rightarrow (S(TM)\otimes E)_{x}$ be parallel
translation about $\nabla^{S(TM)\otimes E}$ along the unique geodesic from $x_0$ to $x$. If $s\in \Gamma (S(TM)\otimes E)$, then we define
$$\widehat{s}_{x_0}(x)=h(x_0,x)\tau(x,x_0)s(x).\eqno(3.4)$$
We write $\widehat{s}_{x_0}(Y)$ instead of $\widehat{s}_{x_0}({\rm exp}_{x_0}Y).$\\
 \indent Let $\theta_X$ be the one-form associated with $X_M$ which is defined
by $\theta_X(Y)=g(X,Y)$ for the vector field $Y$. Let $\nabla^{S,X}$
be the Clifford connection $\nabla^S-\frac{1}{4}\theta_X$ on the
spinors bundle and $\triangle_X$ be the Laplacian on $S(TM)$
associated with $\nabla^{S,X}$. Let
$\mu(X)(\cdot)=\nabla^{TM}_{\cdot}X_M$. Let $U=\{x\in T_{x_0}M|||x||<\varepsilon\}$, where $\varepsilon$ is smaller than the injectivity radius of the manifold $M$ at $x_0$. Define $\alpha:U\times
{\mathfrak{g}}\rightarrow {\mathbb{C}}$ via the formula
$$\alpha_X(x):=-\frac{1}{4}\int_0^1(\iota({\mathcal{R}})\theta_X)(tx)t^{-1}dt,~~\rho(X,x)=e^{\alpha_X(x)},\eqno(3.5)$$
where ${\mathcal{R}}=\sum_{i=1}^nx_i\frac{\partial}{\partial_{x_i}}.$ Then $\rho(X,0)=1$. Recall ([BGV Lemma 8.13]) that the following identity holds
$$H_X=-g^{ij}(x)(\nabla^X_{\partial_i}\nabla^X_{\partial_j}-\sum_k\Gamma^k_{ij}\nabla^X_{\partial_k})+\frac{1}{4}r_M,\eqno(3.6)$$
 where $r_M$ is the scalar curvature and $\Gamma^k_{ij}$ is the connection coefficient of $\nabla^L$.\\
 \indent Let $\tau^X(x_0,x):(S(TM)\otimes E)_{x_0}\rightarrow (S(TM)\otimes E)_{x}$ be parallel
translation about $\nabla^{S\otimes E, X}$ along the unique geodesic from $x_0$ to $x$. If $s\in \Gamma (S(TM)\otimes E)$, then we define
$$\widehat{s}^X_{x_0}(x)=h(x_0,x)\tau^X(x,x_0)s(x).\eqno(3.7)$$
Then $$\widehat{s}^X_{x_0}(x)=\rho\widehat{s}_{x_0}(x).\eqno(3.8)$$\\
where $\rho=\rho(X,x)$ is defined by (3.5).\\

\noindent{\bf Definition 3.3} Let $p\in S\Sigma(E,X)$ and $s\in \Gamma(S(TM)\otimes E)$, then we define
$$\theta(p)(s)(x_0)=\int_{T_{x_0}M\times T_{x_0}^*M}e^{-\sqrt{-1}\left<Y,\xi\right>}\overline{\theta}(p)(x_0,\xi,X)\widehat{s}_{x_0}^X(Y)dYd\xi
.\eqno(3.9)$$\\

\noindent{\bf Remark.} The operator $\theta(p)$ is well-defined since $\sum_{|\alpha|\geq 0}||p_{j,\alpha}(x,\xi)|||X^\alpha|$ and $e^{|\alpha_X(x)|}$ are convergent.
The operator $\theta(p)$ depends on the choice of the cut off function $h$, but the result does not depend on the cut off function for computations of infinitesimal equivariant Chern-Connes characters. We denote by $Op(E,X)$ all such operators with smoothing operators. \\

\noindent{\bf Definition 3.4} Given $s\in \Gamma(S(TM)\otimes E)$, define $\overline{s}_{x_0}^X(x)=h(x_0,x)\tau^X(x_0,x)s(x_0)$
and $\overline{s}_{x_0}(x)=h(x_0,x)\tau(x_0,x)s(x_0)$, then $\overline{s}_{x_0}^X(x)=\rho^{-1}\overline{s}_{x_0}(x)$. Let
$P\in Op(E,X)$ and $s\in \Gamma(S(TM)\otimes E)$. Define $\sigma(P)\in {\rm End}(E)_{x_0}\otimes \Omega(M)\otimes {\mathbb{C}}[\mathfrak{g}^*]$ by
$$\sigma(P)(x_0,\xi,X)=\overline{\sigma}P_y(e^{\sqrt{-1}\left<{\rm exp}^{-1}_{x_0}(y),\xi\right>}\overline{s}_{x_0}^X(y))|_{y=x_0}.\eqno(3.10)$$\\

\noindent{\bf Lemma 3.5} {\it Let $P=\sum_\alpha P_\alpha X^\alpha\in Op(E,X)$. If $\sum_\alpha||P_\alpha||_1|X^\alpha|$ is convergent, then
$\sigma(P)$ is convergent.}\\

\noindent{\bf Proof.} Since $\sum_\alpha||P_\alpha||_1|X^\alpha|$ and $e^{|\alpha_X(x)|}$ are convergent, this comes from the definition 3.4 and $| e^{\sqrt{-1}\left<{\rm exp}^{-1}_{x_0}(y),\xi\right>}|=1$ and $|h(x_0,x)|\leq 1$ and $\tau(x_0,x)$ being an isometry.~~$\Box$\\

\noindent{\bf Lemma 3.6} {\it Let $Y=\sum c_i\partial_i,~Z=\sum d_j\partial_j$ with $c_i,~d_j\in {\mathbb{R}}$. we have}
$$\sigma(\nabla^X_Y )(x,\xi)=\sqrt{-1}\left<Y,\xi\right>_x,\eqno(3.11)$$
$$\sigma(\nabla^X_Y \nabla^X_Z)(x,\xi)=
-\left<Y,\xi\right>\left<Z,\xi\right>+\frac{1}{4}\left<R^L(Y,Z)\partial_k,\partial_l\right>f^k\wedge f^l+\frac{1}{4}\left<\mu^X(Y),Z\right>,\eqno(3.12)$$
{\it where $f^k$ is the dual base of $\partial_k$.}\\

\noindent{\bf Proof.} By the definition 3.4, We have
$$\sigma(\nabla^X_Y )(x_0,\xi)=\overline{\sigma}[\nabla^X_Y(e^{\sqrt{-1}\left<{\rm exp}^{-1}_{x_0}(y),\xi\right>}
\rho^{-1}\overline{s}_{x_0}(y))]|_{y=x_0}.\eqno(3.13)$$
By $$
(d-\frac{1}{4}\theta_X)_{\partial_j}(\rho^{-1})|_{x=x_0}=0;~~\nabla_Y( \overline{s}_{x_0}(x))|_{x=x_0}=0,\eqno(3.14)$$
similarly to the computations of Example 1 in [BF], we get (3.11).
We know that $\rho\nabla^X_Y \nabla^X_Z\rho^{-1}=\rho\nabla^X_Y\rho^{-1}\rho \nabla^X_Z\rho^{-1}$. By the appendix II in [ABP], we have
$$\nabla_Y \nabla_Z\overline{s}_{x_0}(y)|_{y=x_0}=\frac{1}{4}\left<R^L(Y,Z)\partial_k,\partial_l\right>f^k\wedge f^ls(x_0).\eqno(3.15)$$
 In the
trivialization of $S(TM)$, the conjugate
$\rho(X,x)(\nabla^{S,X}_{\partial_i})\rho(X,x)^{-1}$ is given by Lemma 8.13 in [BGV] which is
$$\rho(X,x)(\nabla^{S,X}_{\partial_i})\rho(X,x)^{-1}=\partial_i+
\frac{1}{4}\sum_{j,a<b}\left<R(\partial_i,\partial_j)e_a,e_b \right>
c(e_a)c(e_b)x^j-\frac{1}{4}\mu_{ij}^M(X)x^j$$
$$+\sum_{j<k}f_{ijk}(x)c(e_j)c(e_k)+g_i(x)+\left<h_i(x),X\right>,\eqno(3.16)$$
where $f_{ijk}(x)=O(|x|^2),~g_i(x)=O(|x|)$, and $h_i(x)=O(|x|^2).$
By (3.15) and (3.16), similarly to the computations of Example 2 in [BF], we have (3.12).~~$\Box$\\

\noindent{\bf Proposition 3.7} {\it The following equality holds}
$$\sigma(H_X)=|\xi|^2+\frac{1}{4}r_M.\eqno(3.17)$$
{\it The operator $t^2H_X$ is an asymptotic pseudodifferential operator (see Definition
3.5 in [BF]).}\\

\noindent{\bf Proof.} By Lemma 3.6 and (3.6) and $g^{ij}(x_0)=\delta^{ij}$, $\Gamma^k_{ij}(x_0)=0$ and $R^L(Y,Y)=\left<\mu^X(Y),Y\right>=0$, we get
Proposition 3.7.~~$\Box$\\

\noindent {\bf Definition 3.8} If $p(x,\xi,X)\in S\Sigma(E,X)$, then
$$p_t(x,\xi,X)=\sum^{{\rm dim}M}_{j=0}(\sum_{|\alpha|\geq 0}p_{j,\alpha}(x,t\xi)t^{2|\alpha|}X^\alpha)\otimes \omega_jt^{j},\eqno(3.18)$$\\

 Let $\psi_t: X\rightarrow \frac{X}{t}$ be the rescaling operator on the Lie algebra.\\

\noindent{\bf Theorem 3.9} {\it For $P=\sum P_\alpha X^\alpha\in OP(S\Sigma^{-\infty}(E,X))$ and $t>0$, then}
$$\psi^2_t{\rm Tr_s}(P)=(2\pi)^{-n}(\frac{2}{\sqrt{-1}})^{\frac{n}{2}}\int_M\int_{T^*_{x_0}M}{\rm Tr_s}\sigma( P)_{\frac{1}{t}}(x_0,\xi)d\xi dx.\eqno(3.19)$$
{\it If $P=P_t$ and $P_t $ is an asymptotic pseudodifferential operator and $\sigma(P_t)(x,\xi)$ tends to zero when $|\xi|$ tends to infinity
, then}
$$\psi^2_t{\rm Tr_s}(P_t)=b_0+O(t),\eqno(3.20)$$
{\it where $b_0$ is a constant.}\\

\noindent{\bf Proof.} By Theorem 3.7 in [Ge2], we have for any $s>0$ that
$${\rm Tr_s}(P_\alpha)=(2\pi)^{-n}(\frac{2}{\sqrt{-1}})^{\frac{n}{2}}\int_M\int_{T^*_{x_0}M}{\rm Tr_s}\sigma_G(P_\alpha)_{s}(x_0,\xi)d\xi dx,\eqno(3.21)$$
where
$$\sigma_G(P)(x_0,\xi,X)=\overline{\sigma}P_y(e^{\sqrt{-1}\left<{\rm exp}^{-1}_{x_0}(y),\xi\right>}\overline{s}_{x_0}(y))|_{y=x_0}.\eqno(3.22)$$
Since $\sigma(P_t)(x,\xi)$ tends to zero when $|\xi|$ tends to infinity, by using the equality which will be proved in the following lemma 3.11
$$\int_{T^*_{x_0}M}{\rm Tr_s}\sigma_G(P)_{s}(x_0,\xi)d\xi dx
=\int_{T^*_{x_0}M}{\rm Tr_s}\sigma_G(\rho P\rho^{-1})_{s}(x_0,\xi)d\xi dx,\eqno(3.23)$$
we have for $\rho(x_0)=1$ that
$$\int_{T^*_{x_0}M}{\rm Tr_s}\sigma_G(P_\alpha)_{s}(x_0,\xi)d\xi=\int_{T^*_{x_0}M}{\rm Tr_s}\sigma_G(P_\alpha \rho^{-1})_{s}(x_0,\xi)d\xi$$
$$~~~~~~~~~~~~~~~~~~~~~~~~~~~=\int_{T^*_{x_0}M}{\rm Tr_s}\sigma(P_\alpha)_{s}(x_0,\xi)d\xi.\eqno(3.24)$$
So
$${\rm Tr_s}(P_\alpha X^\alpha)=(2\pi)^{-n}(\frac{2}{\sqrt{-1}})^{\frac{n}{2}}\int_M\int_{T^*_{x_0}M}{\rm Tr_s}\sigma(P_\alpha
\frac{X^\alpha}{s^{2|\alpha|}})_{s}(x_0,\xi)d\xi dx.\eqno(3.25)$$
Let $s=\frac{1}{t}$, then
$${\rm Tr_s}(P_\alpha X^\alpha)=(2\pi)^{-n}(\frac{2}{\sqrt{-1}})^{\frac{n}{2}}\int_M\int_{T^*_{x_0}M}{\rm Tr_s}\sigma(P_\alpha
X^\alpha t^{2|\alpha|})_{\frac{1}{t}}(x_0,\xi)d\xi dx.\eqno(3.26)$$
So
$$\psi^2_t{\rm Tr_s}(P_\alpha X^\alpha)=(2\pi)^{-n}(\frac{2}{\sqrt{-1}})^{\frac{n}{2}}\int_M\int_{T^*_{x_0}M}{\rm Tr_s}\sigma(P_\alpha X^\alpha)_{\frac{1}{t}}(x_0,\xi)d\xi dx.\eqno(3.27)$$
By taking the sum $\sum_\alpha$, we get (3.19). By Definitions 3.8 and Definition 3.5 in [BF], for the  asymptotic pseudodifferential operator  $P_t $, we have
$$\sigma( P_t)=\sum_{l=0}^{+\infty}t^lp_l(x,\xi,X)_t,\eqno(3.28)$$
so
$$\sigma(P_t)_{\frac{1}{t}}=\sum_{l=0}^{+\infty}t^lp_l(x,\xi,X),\eqno(3.29)$$
By (3.19) and (3.29), we get (3.20).~~$\Box$\\

Let $\mu^M$ be the Riemannian moment of $X$ defined by
  $  \mu^M(X)Y=-\nabla_YX^M.$
  Let $F^M_{\mathfrak{g}}(X)=\mu^M+R$ be the equivariant Riemannian curvature of $M$. The equivariant $\widehat{A}$-genus of the tangent bundle of $M$
  is defined by$
\widehat{A}(F^M_{\mathfrak{g}}(X))={\rm det}\left(\frac{F^M_{\mathfrak{g}}(X)/2}{{\rm sinh}(F^M_{\mathfrak{g}}(X)/2
 )}\right)^{\frac{1}{2}}
$.\\

\noindent{\bf Theorem 3.10} {\it When $2k\leq {\rm dim}M$ and $X$ is
small which means that $||X_M||$ is sufficiently small, then for $f^j\in C^{\infty}_G(M)$,}
$${\rm lim}_{t\rightarrow 0}\psi_t{{\bf {\rm ch}}^{2k}}(\sqrt{t}D,X)(f^0,\cdots ,
f^{2k})$$
$$=\frac{1}{(2k)!}{(2\pi\sqrt{-1})}^{-n/2}\int_Mf^0\wedge
df^1\wedge\cdots\wedge
df^{2k}\widehat{A}(F^M_{\mathfrak{g}}(X))d{\rm
Vol}_M.\eqno(3.30)$$\\

\noindent{\bf Proof.} In Theorem 3.9, let $P_t=t^{2k}f^0e^{-\sigma_0t^2H_X}c(df^1)\cdots c(df^{2k})e^{-\sigma_{2k}t^2H_X}$, then by Proposition 3.7,
similarly to Lemma 3.13 in [BF], we have $P_t$ is an asymptotic pseudodifferential operator. By (3.20) and taking the $J$-jet, we have
$${\rm lim}_{t\rightarrow 0}\psi_t^2{\rm Tr_s}(P_t)_J=b_{0,J}.\eqno(3.31)$$
By Theorem 2.12 in [Wa], we have
$${\rm lim}_{t\rightarrow 0}\psi_t^2{\rm Tr_s}(P_t)_J=\frac{1}{(2k)!}{(2\pi\sqrt{-1})}^{-n/2}\int_Mf^0\wedge
df^1\wedge\cdots\wedge
df^{2k}\widehat{A}(F^M_{\mathfrak{g}}(X))_Jd{\rm
Vol}_M.\eqno(3.32)$$
By (3.31) and (3.32) and when $J$ goes to infinity, we obtain
$$b_0=\frac{1}{(2k)!}{(2\pi\sqrt{-1})}^{-n/2}\int_Mf^0\wedge
df^1\wedge\cdots\wedge
df^{2k}\widehat{A}(F^M_{\mathfrak{g}}(X))d{\rm
Vol}_M.\eqno(3.33)$$\\
By (3.20) and (3.33), when $t$ goes to zero, we get (3.30).~~$\Box$\\

\noindent {\bf Lemma 3.11} {\it The equality (3.23) holds.}\\

\noindent{\bf Proof.} Considering the equalities (70) and (71) in [Pf] (Note that these formulas hold for any pseudodifferential operators defined by (3.22) and not only for asymptotic pseudodifferential operators), let $N=0$, then
$$\sigma_G(P\rho^{-1})(x,\xi)=\sigma_G(P)(x,\xi)\rho^{-1}+r_0(\xi)
.\eqno(3.34)$$
where $ r_0(\xi)$ is defined by
$$r_0(\xi)=\frac{\sqrt{-1}}{(2\pi)^n}\sum^n_{j=1}\int_{T^*_{x_0}M(y)}\int^1_0\frac{\partial}{\partial y_j}a(\xi+sy)ds\cdot y_j[{\mathcal{F}}
(f^\psi)](y)dy,\eqno(3.35)$$
and the Fourier transform ${\mathcal{F}}$ and $f^\psi$ are defined by (7) and (8) in [Pf] respectively, $a$ is the symbol of $P$. By (0.2) in [Ge2], we have the leading symbol of $e^{-t^2D^2}$ is $e^{-t^2|\xi|^2}$. As in (2.4), using
the Duhamel principle, we expanse the operator $P_t$ and the leading symbol of $P_t$ is the product of $e^{-t^2|\xi|^2}$ and a polynomial
 on $\xi$. Without loss of generality, we assume $a=e^{-|\xi|^2}$. The following two well-known theorems are necessary:\\

\noindent I. {\it Let $f(x,y)$ be continues on the domain $x\geq a,~y\geq b$ and $\int^{+\infty}_bf(x,y)dy$ be uniformly convergent about
$x$ on any finite interval included in $[a ,+\infty]$ and  $\int^{+\infty}_af(x,y)dx$ be uniformly convergent about
$y$ on any finite interval included in $[b ,+\infty]$. We assume that the integral $\int^{+\infty}_b[\int^{+\infty}_a|f(x,y)|dx]dy$
or $\int^{+\infty}_a[\int^{+\infty}_b|f(x,y)|dy]dx$ exists, then}
$$\int^{+\infty}_a[\int^{+\infty}_bf(x,y)dy]dx=\int^{+\infty}_b[\int^{+\infty}_af(x,y)dx]dy={\rm finite~number}.\eqno(3.36)$$\\

\noindent II. {\it There exists $\beta>0$, such that $|f(x,y)|\leq F(x)$ for any $x>\beta$ and $y\in I$ and that $\int^{+\infty}_aF(x)dx$
exists, then $\int^{+\infty}_af(x,y)dx$ is uniformly convergent.}\\

By (3.35), we consider
$$\int_{T^*_{x_0}M(\xi)}r_0(\xi)d\xi=\frac{\sqrt{-1}}{(2\pi)^n}\sum^n_{j=1}\int_{T^*_{x_0}M(\xi)}\int_{T^*_{x_0}M(y)}
\int^1_0\frac{\partial a}{\partial \xi_j}|_{\xi+sy}sds\cdot y_j[{\mathcal{F}}
(f^\psi)](y)dyd\xi.\eqno(3.37)$$
 Since the Schwartz function $[{\mathcal{F}}
(f^\psi)](y)$ is integral on ${T^*_{x_0}M(y)}$, we take some estimates on the right hand side of (3.37) in the polar coordinates of ${T^*_{x_0}M(\xi)}$ and ${T^*_{x_0}M(y)}$ and then we can verify that the right hand side of (3.37) satisfies the conditions of Theorem I.
 Using
$\int_{T^*_{x_0}M(\xi)}\frac{\partial}{\partial_{\xi_j}}[e^{-|\xi|^2}\xi^\beta] d\xi=0$
and (3.37), we get
 $\int_{T^*_{x_0}M}r_0(\xi)=0.$ Therefore we get (3.23).~~$\Box$\\

 Let $${\rm Ch}({\rm
Im}(p))=\sum_{k=0}^{\infty}(-\frac{1}{2\pi\sqrt{-1}})^k\frac{1}{k!}{\rm
Tr}[p(dp)^{2k}].\eqno (3.38)$$
We have\\

\noindent{\bf Corollary 3.12} {\it When $X$ is
small, then}
$${\rm Ind}_{e^{-X}}(D_{{\rm Im}p,+})
={(2\pi\sqrt{-1})}^{-n/2}\int_M\widehat{A}(F^M_{\mathfrak{g}}(X)){\rm Ch}({\rm Im}p).\eqno(3.39)$$\\

\noindent{\bf Proof.} Using the same discussions as those in [GS], we have the homotopy property of ${\rm ch}^*(D,X)$ for ${t}D_{-X}$. So
by (2.17), we have
$${\rm Ind}_{e^{-X}}(D_{{\rm Im}p,+}) =\left<{\rm ch}^*({t}D_{-X}),{\rm ch}(p)\right>,\eqno(3.40)$$
where $${{\bf {\rm ch}}^{2k}}(tD_{-X})(f^0,\cdots ,
f^{2k}):
=t^{2k}\int_{\triangle_{2k}}{\rm
Str}\left[ e^{-{
\mathfrak{L}}_X}f^0e^{-\sigma_0t^2(D+\frac{1}{4}c(X))^2}c(df^1)\right.$$
$$\left.\cdot e^{-\sigma_1t^2(D+\frac{1}{4}c(X))^2}\cdots c(df^{2k})e^{-\sigma_{2k}t^2(D+\frac{1}{4}c(X))^2}\right]d{\rm Vol}_{\Delta_{2k}},\eqno(3.41)$$
In (3.40), let $e^{-X}=e^{-t^2X}$ and use $(\psi_t)^2$ acting on (3.40), then we get
$${\rm Ind}_{e^{-X}}(D_{{\rm Im}p,+}) =\left<{\widetilde{{\rm ch}}}^*(tD,{X}),{\rm ch}(p)\right>,\eqno(3.42)$$
where$${\widetilde{{\rm ch}}}^{2k}(tD,{X})(f^0,\cdots ,
f^{2k}):
=t^{2k}\int_{\triangle_{2k}}{\rm
Str}\left[ f^0e^{-\sigma_0t^2H_{\frac{X}{t^2}}}c(df^1)\right.$$
$$\left.\cdots c(df^{2k})e^{-\sigma_{2k}t^2H_{\frac{X}{t^2}}}\right]d{\rm Vol}_{\Delta_{2k}}.\eqno(3.43)$$
 Since ${\rm Ind}_{e^{-X}}(D_{{\rm Im}p,+}) $ is independent of $t$,
taking the limit as $t\rightarrow 0$ in (3.42), we get by Theorem 3.10 that
$${\rm Ind}_{e^{-X}}(D_{{\rm Im}p,+})={(2\pi\sqrt{-1})}^{-n/2}\int_M\widehat{A}(F^M_{\mathfrak{g}}(X)){\rm Ch}({\rm Im}p)d{\rm
Vol}_M.\eqno(3.44)$$
~~$\Box$\\

\section{The infinitesimal equivariant eta cochains}

\quad In this Section, we prove the limit of truncated infinitesimal equivariant eta cochains exists when $J$ goes to infinity.
By the Duhamel principle and (2.5), we have
\begin{eqnarray*}
&&||D_{-X}e^{-utH_X}||_{u^{-1}}\\
&\leq&\sum_{m\geq 0}(ut)^m\int_{\triangle_m}||D_{-X}(1+D^2)^{-\frac{1}{2}}|| ||(1+D^2)^{\frac{1}{2}}e^{-\frac{\sigma_0}{2}utD^2}||
||e^{-\frac{\sigma_0}{2}utD^2}||_{(u\sigma_0)^{-1}}\\
&&\cdot||F_X(1+D^2)^{-\frac{1}{2}}||||(1+D^2)^{\frac{1}{2}}e^{-\frac{\sigma_1}{2}utD^2}||
||e^{-\frac{\sigma_1}{2}utD^2}||_{(u\sigma_1)^{-1}}\\
&&\cdots
||F_X(1+D^2)^{-\frac{1}{2}}||||(1+D^2)^{\frac{1}{2}}e^{-\frac{\sigma_m}{2}utD^2}||
||e^{-\frac{\sigma_m}{2}utD^2}||_{(u\sigma_m)^{-1}}d\sigma\\
&\leq &||D_{-X}(1+D^2)^{-\frac{1}{2}}||(eut)^{-\frac{1}{2}}\sum_{m\geq 0}(e^{-1}ut||F_X(1+D^2)^{-\frac{1}{2}}||^2)^{\frac{m}{2}}\\
&&\cdot
e^{\frac{ut}{2}}({\rm tr}e^{-\frac{t}{2}D^2})^u\int_{\triangle_m}\sigma^{-\frac{1}{2}}_0\cdots \sigma^{-\frac{1}{2}}_md\sigma\\
&\leq&||D_{-X}(1+D^2)^{-\frac{1}{2}}||(ut)^{-\frac{1}{2}}2e^{\frac{ut}{2}}\{1+[||F_X(1+D^2)^{-\frac{1}{2}}||^{2}e^{-1}\pi ut]^{\frac{1}{2}}\}\\
&&e^{||F_X(1+D^2)^{-\frac{1}{2}}||^{2}\pi ut}
({\rm tr}e^{-\frac{t}{2}D^2})^u,
~~~~~~~~~~~~~~~~~~~~~~~~~~~~~~~~~~~~~~~~~~~~~~~~~~~~~~~~~(4.1)
\end{eqnarray*}
where
$$\int_{\triangle_m}\sigma_0^{-\frac{1}{2}}\cdots \sigma_{m}^{-\frac{1}{2}}d\sigma
=\frac{\pi^{\frac{m+1}{2}}}{\Gamma(\frac{m}{2}+1)},\eqno(4.2)$$
and
$$\frac{\pi^{\frac{m+1}{2}}}{\Gamma(\frac{m}{2}+1)}=\frac{\pi^{\frac{m+1}{2}}}{(\frac{m}{2})!},~~ {\rm when~} m~ {\rm is~ even},\eqno(4.3)$$
$$\frac{\pi^{\frac{m+1}{2}}}{\Gamma(\frac{m}{2}+1)}\leq\frac{2\pi^{\frac{m}{2}}}{(\frac{m-1}{2})!},~~{\rm when~} m~ {\rm is~ odd}.\eqno(4.4)$$
\indent Now let $M$ be a compact oriented odd dimensional Riemannian
manifold without boundary with a fixed spin structure and $S$ be the
bundle of spinors on $M$. The fundamental setup consists with that on page 2. Let $K_t=\sqrt{t}(D+\frac{c(X)}{4t})$, then $\frac{dK_t}{dt}=\frac{1}{2\sqrt{t}}
D_{\frac{X}{t}}.$ For $a_0,\cdots,a_{2k}\in C_G^{\infty}(M)$, we define
the infinitesimal equivariant
 cochain ${\rm ch}^{2k}_X({{K}}_t,\frac{d{{K}}_t}{dt})$ by the formula:\\
$${\rm ch}^{2k}_X({{K}}_t,\frac{d{{K}}_t}{dt})(a_0,\cdots ,
a_{2k})~~~~~~~~~~~~~~~~~~~~~~$$
$$=
\sum^{2k}_{j=0}(-1)^j\langle a_0,[{{K}}_t,a_1],\cdots ,[{{K}}_t,a_{j}],\frac{d{{K}}_t}{dt},[{{K}}_t,a_{j+1}]
,\cdots , [{{K}}_t,a_{2k}]\rangle _t(X).\eqno(4.5)$$
 If $A_j~ (0\leq j\leq q)$ are operators on $\Gamma (M,S(TM))$, we define\\
\bigskip
$$\langle A_0,\cdots , A_q\rangle_t(X)=\int_{\triangle_q}{\rm tr}[e^{-L_X}A_0e^{-\sigma_0K^2_t}A_1e^{-\sigma_1K^2_t}\cdots A_q
e^{-\sigma_qK^2_t}]d\sigma,\eqno(4.6)$$ where
$\triangle_q=\{(\sigma_0,\cdots,\sigma_q)|\sigma_0+\cdots+\sigma_q=1,~\sigma_j\geq 0\}$ is a simplex in ${\bf R^q}$
and $L_X$ is the Lie derivative generated by $X$ on the spinors bundle.
\\
\indent Formally, {\bf the infinitesimal equivariant eta cochain} for the odd dimensional manifold is defined to be an even cochain sequence
by the formula:\\
$${\eta}^{2k}_X({{D}})=\frac{1}{\sqrt{\pi}}\int^{\infty}_{0}{\rm ch}^{2k}_X
({{K}}_t,\frac{d{{K}}_t}{dt})dt,
\eqno(4.7)$$
 Then
${\eta}^0_{X}(D)(1)$ is the half of the
infinitesimal equivariant eta invariant defined by Goette in [Go].
In order to prove that the above expression is well defined, it is necessary to check the
integrality near the two ends of the integration.
 Firstly, the regularity at infinity comes from the
following lemma.\\

\noindent {\bf Lemma 4.1}~~{\it For
$a_0,\cdots,a_{2k}\in C^{\infty}_G(M)$, we have}
$${\rm ch}^{2k}_X({{K}}_t,\frac{d{{K}}_t}{dt})(a_0,\cdots ,
a_{2k})=O(t^{-\frac{3}{2}}),~~{\rm as}~t\rightarrow\infty.\eqno(4.8)$$

\noindent{\bf Proof.} Let $L_0$ be a fixed large number. Then
$\frac{1}{\Gamma(\frac{1}{2})}\int^{L_0}_{\varepsilon}{\rm ch}^{2k}_X({{K}}_t,\frac{d{{K}}_t}{dt})(a_0,\cdots ,
a_{2k})dt$ is well-defined
 by Lemma 2.2 and (4.1). Similarly to Lemma 2.2 and (4.1), we know that Lemma 3.5 in [Wa] holds when $J$ goes to infinity. So $\frac{1}{\Gamma(\frac{1}{2})}\int^{\infty}_{L_0}{\rm ch}^{2k}_X({{K}}_t,\frac{d{{K}}_t}{dt})dt$ is well-defined and Lemma 4.1
 holds.~~~~~~~$\Box$\\

\indent Next, we prove the regularity at zero. Let ${ F_*}=D^2_{-X}$ and $\widehat{{ F_*}}=H_X-dtD_X$ where $dt$ is an
auxiliary Grassmann variable as shown in [BiF]. Then $t\psi_t\widehat{{F_*}}=tH_{\frac{X}{t}}-2t^{\frac{3}{2}}dt\frac{dK_t}{dt}.$
Let
$${\rm ch}^{2k}(\widehat{{F_*}})(a_0,\cdots ,
a_{2k})
=t^k\int_{\triangle_{2k}}\psi_t{\rm tr}[a_0e^{-t\sigma_0\widehat{{ F_*}}}[D,a_1]\cdots [D,a_{2k}]
e^{-t\sigma_{2k}\widehat{{ F_*}}}]d\sigma,\eqno(4.9)$$
$${\rm ch}^{2k}({{ F_*}})(a_0,\cdots ,
a_{2k})
=t^k\int_{\triangle_{2k}}\psi_t{\rm tr}[a_0e^{-t\sigma_0H_X}[D,a_1]\cdots [D,a_{2k}]
e^{-t\sigma_{2k}H_X}]d\sigma.\eqno(4.10)$$
By the Duhamel principle and $dt^2=0$, we have
$$e^{-t\sigma_j\psi_t\widehat{ F_*}}=e^{-t\sigma_j\psi_tH_X}+\int^1_0e^{-(1-a)t\sigma_j\psi_tH_X}(2t^{\frac{3}{2}}dt\frac{dK_t}{dt})
e^{-at\sigma_j\psi_tH_X}d(\sigma_ja)$$
$$=e^{-t\sigma_j\psi_tH_X}+2t^{\frac{3}{2}}dt\int^{\sigma_j}_0e^{-(\sigma_j-\xi)t\psi_tH_X}\frac{dK_t}{dt}e^{-t\xi\psi_tH_X}d\xi
\eqno(4.11)$$
By (4.5) and (4.9)-(4.11) and $dt^2=0$, we get
$${\rm ch}^{2k}(\widehat{{ F_*}})(a_0,\cdots ,
a_{2k})={\rm ch}^{2k}({{F_*}})(a_0,\cdots ,
a_{2k})-2t^{\frac{3}{2}}{\rm ch}^{2k}_X({{K}}_t,\frac{d{{K}}_t}{dt})(a_0,\cdots ,
a_{2k})dt.\eqno(4.12)$$\\

\noindent {\bf Lemma 4.2} {\it The following estimate holds}\\
$${\rm ch}^{2k}_X({{K}}_t,\frac{d{{K}}_t}{dt})\sim O(1) ~~{\rm when}~t\rightarrow 0.\eqno(4.13)$$\\

\noindent{\bf Proof.} By (4.12), we only need to prove
$${\rm ch}^{2k}(\widehat{{ F_*}})(a_0,\cdots ,
a_{2k})-{\rm ch}^{2k}({{F_*}})(a_0,\cdots ,
a_{2k})=O(t^{\frac{3}{2}})dt.\eqno(4.14)$$
Let
$$Q_{\widehat{ F_*}}=a_0(\widehat{ F_*}+\partial_t)^{-1}c(da_1)\cdots c(da_{2q})(\widehat{ F_*}+\partial_t)^{-1},\eqno(4.15)$$
$$
Q_{{ F_*}}=a_0({ F_*}+\partial_t)^{-1}c(da_1)\cdots c(da_{2q})({ F_*}+\partial_t)^{-1}.\eqno(4.16)$$
By using Lemma 8.4 in [PW], we have
$$t^q\psi_t[a_0e^{-t\sigma_0\widehat{{ F_*}}}[D,a_1]\cdots [D,a_{2q}]
e^{-t\sigma_{2q}\widehat{{ F_*}}}](x,y)=t^{-q}\psi_tK_{Q_{\widehat{ F_*}}}(x,y,t);\eqno(4.17)$$
$$t^q\psi_t[a_0e^{-t\sigma_0H_X}[D,a_1]\cdots [D,a_{2q}]
e^{-t\sigma_{2q}H_X}](x,y)=t^{-q}\psi_tK_{Q_{ F_*}}(x,y,t).\eqno(4.18)$$
So we only need to prove
$$t^{-q}\psi_t{\rm tr}\left[K_{Q_{\widehat{ F_*}}}(x,x,t)-K_{Q_{ F_*}}(x,x,t)\right]=O(t^{\frac{3}{2}})dt.\eqno(4.19)$$
By the trace property, we have
$$t^{-q}\psi_t{\rm tr}\left[K_{Q_{\widehat{ F_*}}}(x,x,t)-K_{Q_{ F_*}}(x,x,t)\right]$$
$$=t^{-q}\psi_t{\rm tr}\left[K_{Q_{h\rho \widehat{F_*}(h\rho)^{-1}}}(x,x,t)-K_{Q_{\rho H_X\rho^{-1}}}(x,x,t)\right].\eqno(4.20)$$
By (3.15), (3.18) and (3.24) in [Wa] and $dt^2=0$ where we use $dt$ instead of $z$ in [Wa], we have
$$t^{-q}\psi_t\left[{Q_{h\rho \widehat{F_*}(h\rho)^{-1}}}-{Q_{\rho H_X\rho^{-1}}}\right]$$
$$=-t^{-q}dt\psi_t\left[a_0(\partial_t+\rho H_X\rho^{-1})^{-1}u(\partial_t+\rho H_X\rho^{-1})^{-1}\right.$$
$$\cdot c(da_1)\cdots c(da_{2q})(\partial_t+\rho H_X\rho^{-1})^{-1}$$
$$+\cdots +
a_0(\partial_t+\rho H_X\rho^{-1})^{-1}\cdots
c(da_{2q})$$
$$\left.\cdot(\partial_t+\rho H_X\rho^{-1})^{-1}u(\partial_t+\rho H_X\rho^{-1})^{-1}\right]
.\eqno(4.21)$$
By $O_G(u)\leq 0$ and $O_G((\partial_t+\rho H_X\rho^{-1})^{-1})=-2$, we have
$$O_G\left[(\partial_t+\rho H_X\rho^{-1})^{-1}u(\partial_t+\rho H_X\rho^{-1})^{-1}\right.$$
$$\left.\cdot c(da_1)\cdots c(da_{2q})(\partial_t+\rho H_X\rho^{-1})^{-1}\right]=-2q-4,\eqno(4.22)$$
which has odd Clifford elements. When we drop off the truncated order $J$ in Lemma 2.9 in [Wa] and consider the convergent series on $X$
as in the definition 3.2, we know that Lemma 2.9 in [Wa] holds for our operator in (4.22). By (4.20)-(4.22) and Lemma 2.9 1) in [Wa] for $j=n$ and $m=-2q-4$, we get (4.19).~~$\Box$\\

\noindent{\bf Remark.} Similarly to Proposition 1.2 in [Wu], We use the symbol calculus about the connection $\nabla^X$ in Section 3 instead of the Getzler symbol calculus in Proposition 1.2 in [Wu], then we can give another proof of Lemma 4.2.\\

\indent Again Proposition 3.8 in [Wa] holds, we have\\

 \noindent {\bf Proposition 4.3} {\it Assume that $D$ is invertible
with $\lambda$ being the smallest positive eigenvalue of $|D|$ and $||dp||<\lambda,$
 then
the pairing $\langle\eta^*_{X}(D),{\rm ch}_*(p)\rangle$ is well-defined.}\\

\indent We also have the following theorem.\\

\noindent {\bf Theorem 4.4} {\it Assume $D$ is invertible and $||dp||<\lambda$ where $\lambda$ is the smallest eigenvalue of $|D|$, then we have}
 $$\frac{1}{2}\eta_X(p(D\otimes
 I_r)p)=\langle\eta^*_X(D),{\rm ch}_*(p)\rangle
 ,\eqno(4.23)$$
{\it where $\eta_X(p(D\otimes
 I_r)p)$ is the Goette's infinitesimal equivariant eta invariant.}\\

 \noindent{\bf Proof.} We still use the same notations and discussions after Proposition 3.8 in [Wa]. The difference is that we add $\psi_t$ in the definition of $A$. That is,
let $A=d_{(u,s,t)}+\psi_t\widetilde{{\bf D_{-X}}}$ be a
superconnection on the trivial infinite dimensional superbundle with the
base $[0,1]\times {\bf R}\times (0,+\infty)$ and the fibre $H\otimes {\bf
C^r}\oplus H\otimes {\bf C^r}.$ Then we have
$$A^2=t\psi_t{\bf
D}_{-X,u}^2-s^2/4-(1-u)t^{\frac{1}{2}}s\sigma[{\bf D},p]+ds\sigma
(p-\frac{1}{2})+t^{\frac{1}{2}}du(2p-1)[{\bf
D},p]+\frac{dt}{2t^{\frac{1}{2}}}\psi_t{\bf
D}_{X,u}.\eqno(4.24)$$
Since we prove the regularity at zero, we can take $\varepsilon=0$ in (3.41)-(3.45) in [Wa]. By
the following lemma, Theorem 4.4 can be proved.~~$\Box$\\

\noindent{\bf Lemma 4.5}~{\it Let $D_u=D+u(2p-1)[D,p]$ for $u\in [0,1]$. We assume that $D$ be invertible and $||dp||<\lambda$, then we have}
${\eta}_X(D_0)={\eta}(D_1)$.\\

\noindent{\bf Proof.} By $||dp||<\lambda$, then $D_u=D+u(2p-1)[D,p]$ is invertible for $u\in [0,1]$. Similar to the discussions of Proposition 4.4 in [Wu], the infinitesimal equivariant eta invariant of $D_u$ is well defined. So ${\eta}_X(D_u)$ is smooth. Let $A=(2p-1)dp$. Then by the definition of the infinitesimal equivariant eta invariant and the Duhamel principle, we have
$$\frac{d}{d u}{\eta}_X(D_u)=\frac{1}{\sqrt{\pi}}\int^{+\infty}_0{\rm tr}[e^{-X}Ae^{-tD^2_{-\frac{X}{t},u}}]d\sqrt{t}+L,\eqno(4.25)$$
where $$
L=-\frac{t^{\frac{1}{2}}}{2\sqrt{\pi}}\int^{+\infty}_0\int^1_0{\rm tr}\left\{e^{-X}D_{\frac{X}{t},u}e^{-(1-s)tD^2_{-\frac{X}{t},u}}[D_{-\frac{X}{t},u},A]_+
e^{-stD^2_{-\frac{X}{t},u}}ds\right\}dt.\eqno(4.26)$$
By the trace property and direct computations, then
$$\frac{\partial}{\partial t}(\sqrt{t}D_u+\frac{c(X)}{4\sqrt{t}})^2=\frac{1}{2}[D_u+\frac{c(X)}{4t},D_u-\frac{c(X)}{4t}]_+,\eqno(4.27)$$
$$\int^1_0{\rm tr}\left\{Ae^{-(1-s)tD^2_{-\frac{X}{t},u}}[D_{-\frac{X}{t},u},D_{\frac{X}{t},u}]_+
e^{-stD^2_{-\frac{X}{t},u}}\right\}ds$$
$$=\int^1_0{\rm tr}\left\{D_{\frac{X}{t},u}e^{-(1-s)tD^2_{-\frac{X}{t},u}}[D_{-\frac{X}{t},u},A]_+
e^{-stD^2_{-\frac{X}{t},u}}\right\}ds.\eqno(4.28)$$
By using the Duhamel principle and the Leibniz rule and (4.26)-(4.28), we get
$$\frac{\partial}{\partial u}\psi_t{\rm tr}[D_{X,u}e^{-t(D^2_{-X,u}+L_X)}]d\sqrt{t}=
\frac{\partial}{\partial t}{\rm tr}[t^{\frac{1}{2}}e^{-X}Ae^{-t{D^2_{-\frac{X}{t},u}}}]dt.\eqno(4.29)$$
So
$$\frac{d}{d u}{\eta}_X(D_u)=\frac{1}{\sqrt{\pi}}{\rm tr}[t^{\frac{1}{2}}e^{-X}Ae^{-t{D^2_{-\frac{X}{t},u}}}]\left|^{+\infty}_{t=0}\right..\eqno(4.30)$$
As $D_u$ is invertible, then
$${\rm lim}_{t\rightarrow +\infty}{\rm tr}[t^{\frac{1}{2}}e^{-X}Ae^{-t{D^2_{-\frac{X}{t},u}}}]=0.\eqno(4.31)$$
Using Lemma 2.9 in [Wa] for $j=n$ and $m=-1$, similar to the discussions on Line 14 in [Wu, P.164], we have
$${\rm lim}_{t\rightarrow 0}{\rm tr}[t^{\frac{1}{2}}e^{-X}Ae^{-t{D^2_{-\frac{X}{t},u}}}]~~~~~~~~~~~~~~~~~~~~~~~~~~~~~~$$
$$=c_0\int_M\widehat{A}(F^M_{\mathfrak{g}}(X)){\rm tr}\left\{(2p-1)(dp){\rm exp}[\frac{\sqrt{-1}}{2\pi}(A'\wedge A'+dA')]\right\}=0,\eqno(4.32)$$
 where $A'=u(2p-1)dp$. Then by (4.30)-(4.32), Lemma 4.5 is proved.~~$\Box$\\

 \indent Let $N$ be an even-dimensional compact manifold with the boundary $M$. We endow $N$ with a metric which is a product in
a collar neighborhood of $M$. Denote by $D~(D_M)$ the Dirac operator on $N~(M)$. Let $C^{\infty}_*(N)=\{f\in C^{\infty}(N)|f$ is
independent of the normal coordinate $x_n$ near the boundary
$\}.$\\

 \noindent {\bf Definition
4.6}~The infinitesimal equivariant Chern-Connes character on $N$, $\tau_X=\{
\tau^0_X,\tau^2_X,\cdots, \tau^{2q}_X\cdots \}$ is defined by
$$ \tau^{2q}_X(f^0,f^1,\cdot,f^{2q}):= -{\eta}^{2q}_X(D_M)
(f^0|_M,f^1|_M,\cdot,f^{2q}|_M)$$
$$+\frac{1}{(2q)!(2\pi\sqrt{-1})^q}
\int_M \widehat{A}(F^M_{\mathfrak{g}}(X))f^0 df^1\wedge \cdots\wedge
df^{2q},\eqno(4.33)$$
where $f^0,f^1,\cdot,f^{2q}\in C^{\infty}_*(N)$.\\

\indent Similarly to Proposition 4.2 in [Wa1], we have\\

 \noindent {\bf Proposition 4.7}~{\it The infinitesimal equivariant Chern-Connes character is
$b-B$ closed (for the definitions of $b,~B$, see [FGV]). That is, we have}
$$b\tau^{2q-2}_X+B\tau^{2q}_X=0.\eqno(4.34) $$\\

\indent By Proposition 4.3, we have\\

\noindent {\bf Proposition 4.8} {\it Suppose that $D_{M}$ is invertible
with $\lambda$ being the smallest positive eigenvalue of $|D_{M}|$. We assume that
$||d(p|_M)||<\lambda$, then
the pairing $\langle\tau^*_X,{\rm ch}_*(p)\rangle$ is well-defined.}\\

\indent We let $C_1(M)=M\times (0,1],~\widetilde{N}=N\cup_{M\times \{1\}}C_1(M)$
and $\cal{U}$ be a collar neighborhood of $M$ in $N$. For
$\varepsilon>0$, we take a metric $g^{\varepsilon}$ of $\widetilde{N}$ such that
on ${\cal{U}}\cup_{M\times \{1\}}C_1(M)$
$$g^{\varepsilon}=\frac{dr^2}{\varepsilon}+r^2g^{M}.$$
\noindent Let $S=S^+\oplus S^-$ be spinors bundle associated to
$(\widetilde{N},g^{\varepsilon})$ and $H^{\infty}$ be the set
$\{\xi\in\Gamma(\widetilde{N},S)|~\xi ~{\rm and~ its~ derivatives~ are~ zero~
near~ the~ vertex~ of~ cone~}\}.$ Denote by $L^2_c(\widetilde{N},S)$ the
$L^2-$completion of $H^{\infty}$ (similarly define $L^2_{c}(\widetilde{N},S^+)$
and $L^2_{c}(\widetilde{N},S^-)$). Let
$$D_{\varepsilon}:~H^{\infty}\rightarrow H^{\infty};
~~D_{+,\varepsilon}:~H_+^{\infty}\rightarrow H_-^{\infty},$$ be the
Dirac operators associated with $(\widetilde{N},g^{\varepsilon})$ which are
Fredholm operators for the sufficiently small $\varepsilon$. By $||d(p|_M)||<\lambda$, then $pD_Mp$ is invertible. Recall the Goette's infinitesimal equivariant index theorem for the twisting bundle ${\rm Im}p$ with the connection $pd$ in [Go] that
$${\rm Ind}_{e^{-X}}(pD_{+,\varepsilon}p)=\sum_{r=0}^{\infty} \frac{(-1)^r}{r!(2\pi\sqrt{-1})^r}
\int_{N}
\widehat{A}(F^N_{\mathfrak{g}}(X)){\rm
Tr}[p(dp)^{2r}]-\frac{1}{2}{\eta}_X(pD_Mp).\eqno(4.35)$$
By the Stokes theorem and the trace property and $p(dp)^2=(dp)^2p$, we have
$$\int_M\widehat{A}(F^M_{\mathfrak{g}}(X)){\rm tr}[p_M(d_Mp_M)^{2k-1}]=0.\eqno(4.36)$$
By $L_X(p)=\iota_Xd(p)=0$, then $\iota_X[p(dp)^{2k-1}]=0$. By the Stokes theorem and (4.36), we get
$$\int_N\widehat{A}(F^N_{\mathfrak{g}}(X)){\rm tr}[(d_Np_N)^{2k}]=\int_N(d+\iota_X)\left[\widehat{A}(F^N_{\mathfrak{g}}(X)){\rm tr}[p(d_Np_N)^{2k-1}]\right]$$
$$=\int_M\widehat{A}(F^M_{\mathfrak{g}}(X)){\rm tr}[p_M(d_Mp_M)^{2k-1}]=0.\eqno(4.37)$$
By Theorem 4.4 and Definition 4.6 and (2.14) and (4.37), we get\\

\noindent {\bf Theorem 4.9}~{\it Suppose that $D_{M}$ is invertible
with $\lambda$ being the smallest positive eigenvalue of $|D_{M}|$. We assume that
$||d(p|_M)||<\lambda$ and $p\in M_{r\times r}(C^{\infty}_*(N))$, then }
$${\rm Ind}_{e^{-X}}(pD_{+,\varepsilon}p)=\langle\tau^*_X(D),{\rm
ch}_*(p)\rangle.\eqno(4.38)$$\\

\section{ The infinitesimal equivariant Chern-Connes character for a family of Dirac operators}

\quad In this Section, we extend Sections 2, 3 to the family case. Let us recall the definition of the equivariant family Bismut Laplacian.
Let $M$ be a $n+{q}$ dimensional compact connected manifold
 and $B_0$ be a ${q}$ dimensional compact connected
 manifold. Assume that $\pi :M\rightarrow B_0$ is a fibration and $M$ and $B_0$ are oriented. Taking the orthogonal bundle of the vertical bundle $TZ$
 in $TM$ with respect to any Riemannian metric, determines a smooth
 horizontal subbundle $T^HM$, i.e. $TM=T^HM\oplus TZ$.  Recall that $B_0$ is Riemannian, so we can lift the
 Euclidean scalar product $g_{B_0}$ of $TB_0$ to $T^HM$.
And we assume that $TZ$ is endowed with a scalar product $g_Z$. Thus
we can introduce a new scalar product $g_{B_0}\oplus g_Z$ in $TM$.
Denote by $\nabla^L$ the Levi-Civita connection on $TM$ with respect
to this metric. Let $\nabla^{B_0}$ denote the Levi-Civita connection on
$TB_0$ and still denote by $\nabla^{B_0}$ the pullback connection on
$T^HM$. Let $\nabla^Z=P_Z(\nabla^L)$, where $P_Z$ denotes the
projection to $TZ$. Let $\nabla^{\oplus}=\nabla^{B_0}\oplus \nabla^Z$
and $\omega=\nabla^L-\nabla^{\oplus}$ and $T$ be the torsion tensor
of $\nabla^{\oplus}$. Now we assume that the bundle $TZ$ is spin.
Let $S(TZ)$ be the associated spinors bundle and $\nabla^Z$ can be
lifted to give a connection on $S(TZ)$. Let $D$ be the tangent Dirac
operator.\\
\indent Let $G$ be a compact Lie group which acts fiberwise on $M$.
We will consider that $G$ acts as identity on $B_0$. We assume that the action of G lifts to $S(TZ)$
and the $G$-action commutes with $D$. Let $E$ be the vector
bundle $\pi^*(\wedge T^*{B_0})\otimes S(TZ)$. This bundle carries a
natural action $m_0$ of the degenerate Clifford module $Cl_0(M)$.
Define the connection for $X\in {\mathfrak{g}}$ whose Killing vector
field is in $TZ$,
$$\nabla^{E,-X,\oplus}:=\pi^*\nabla^{B_0}\otimes 1+1\otimes \nabla^{S,-X},\eqno(5.1)$$
$$\omega(Y)(U,V):=g(\nabla^L_YU,V)-g(\nabla^{\oplus}_YU,V),\eqno(5.2)$$
$$\nabla^{E,-X,0}_Y:=\nabla^{E,-X,\oplus}_Y+\frac{1}{2}m_0(\omega(Y)),\eqno(5.3)$$
for $Y,U,V\in TM$. Then the equivariant Bismut superconnection acting on $\Gamma(M,\pi^*\wedge(T^*{B_0})\otimes S(TZ))$ is defined by
$${{B}}^{-X}=\sum_{i=1}^nc(e_i^*)\nabla_{e_i}^{E,-X,0}+\sum_{j=1}^{{q}}f_j^*\wedge\nabla_{f_j}^{E,-X,0};~
~B^{-X}=B+\frac{1}{4}c(X).\eqno(5.4)$$
where $e_1,\cdots,e_n$ and $f_1,\cdots,f_q$ are orthonormal basis of $TZ$ and $TB_0$ respectively, and $B$ is the Bismut superconnection defined by
$$\nabla^{E,\oplus}:=\pi^*\nabla^{B_0}\otimes 1+1\otimes \nabla^{S};\eqno(5.5)$$
$$\nabla^{E,0}_Y:=\nabla^{E,\oplus}_Y+\frac{1}{2}m_0(\omega(Y));\eqno(5.6)$$
$${{B}}=\sum_{i=1}^nc(e_i^*)\nabla_{e_i}^{E,0}+\sum_{j=1}^{\overline{q}}c(f_j^*)\nabla_{f_j}^{E,0}.\eqno(5.7)$$
Define the equivariant family Bismut Laplacain as follows:
$$H_{B,X}=({{B}}^{-X})^2+{{L}}^E_X,\eqno(5.8)$$
where ${L}^E_X$ is the Lie derivative induced by $X$ on the bundle $E$.
Then $$H_{B,X}=D^2+F_++{\widetilde{F}}_+,\eqno(5.9)$$
where $D^2_{-X}=D^2+F_+$ and ${\widetilde{F}}_+=H_{B,X}-D^2_{-X}$ is a first order differential operator along the fibre with coefficients in $\Omega_{\geq 1}(B_0)$.\\

\noindent {\bf Definition 5.1} ~The infinitesimal equivariant family JLO cochain ${{\bf
{\rm ch}}}^{
2k}(B,X)$ can be defined by the formula for $f^0,\cdots,f^{2k}$ in $C^{\infty}_G(M)$:\\
$${{\bf {\rm ch}}^{2k}}(B,X)(f^0,\cdots ,
f^{2k}):
=\int_{\triangle_{2k}}{\rm
Str}\left[ f^0e^{-\sigma_0H_{B,X}}c(df^1)e^{-\sigma_1H_{B,X}}\right.$$
$$\left.\cdots c(df^{2k})e^{-\sigma_{2k}H_{B,X}}\right]d{\rm Vol}_{\Delta_{2k}},\eqno(5.10)$$
where Str is taking the trace along the fibre.\\

\indent Similarly to Section 2, we can prove that (5.10) is well-defined and $\left<{\rm ch}^*(B,X), {\rm ch}p\right>$ is convergent by the following lemma.\\

\noindent {\bf Lemma 5.2} ~{\it For any $1\geq u> 0$, we have:}\\
$$||e^{-uH_{B,X}}||_{u^{-1}}\leq C_0e^{||F_X(1+D^2)^{-\frac{1}{2}}||\pi u}
 ({\rm
tr}[e^{-\frac{D^2}{2}}])^u,\eqno(5.11)$$
{\it where the constant $C_0$ is independent of $u$.}\\

\noindent{\bf Proof.} By (5.10) and the Duhamel principle, we have
$$e^{-uH_{B,X}}=e^{-uH_X}+\sum^{{\rm dim}B_0}_{r>0}I_r,\eqno(5.12)$$
where
$$I_r=\int_{\triangle_r}e^{-s_0uH_X}{\widetilde{F}}_+e^{-s_1uH_X}\cdots{\widetilde{F}}_+e^{-s_ruH_X}ds.\eqno(5.13)$$
In (4.1), we use ${\widetilde{F}}_+$ and $su$ instead of $D_{-X}$ and $u$ respectively and let $t=1$, then we have
$$||{\widetilde{F}}_+e^{-suH_X}||_{(su)^{-1}}\leq 2(su)^{-\frac{1}{2}}
||{\widetilde{F}}_+(1+D^2)^{-\frac{1}{2}}||
e^{\frac{su}{2}}$$
$$\cdot\{1+[||F_X(1+D^2)^{-\frac{1}{2}}||^{2}e^{-1}\pi su]^{\frac{1}{2}}\}\\
e^{||F_X(1+D^2)^{-\frac{1}{2}}||^{2}\pi su}
({\rm tr}e^{-\frac{1}{2}D^2})^{su}.\eqno(5.14)$$
By Lemma 2.2 and (5.12)-(5.14) and the H\"{o}lder
inequality, we get Lemma 5.2.~~~~~~$\Box$\\

\indent Similarly to Propositions 4.11 and 4.12 in [BC], we have\\

\noindent {\bf Proposition 5.3} (1) {\it  The infinitesimal equivariant family Chern-Connes character is closed:}
$$(B+b+d_{B_0})({\rm ch}^*(B,X))=0.\eqno(5.15)$$
(2)~ {\it Let $B_{\tau}=B^{-X}+\tau V$ and $V$ is a bounded operator which commutes with $e^{-X}$, then there exists a cochain ${\rm ch}^*(B_{\tau},X,V)$ such that}
$$\frac{d}{d\tau}{\rm ch}^*(B_{\tau},X)=-[b+B+d_{B_0}]{\rm ch}^*(B_{\tau},X,V).\eqno(5.16)$$

\indent By taking $V=(2p-1)[B,p]$, we get\\

\noindent{\bf Theorem 5.4} {\it The following index formula holds in the cohomology of $B_0$}
$${\rm Ch}_{e^{-X}}[{\rm Ind}(D_{{\rm Im}p,+,z})] =\left<{\rm ch}^*(B,X),{\rm ch}(p)\right>.\eqno(5.17)$$

\indent Let $\phi_t$ be the rescaling operator on $\Omega(B_0)$ defined by $dy_j\rightarrow \frac{dy_j}{\sqrt{t}}$ for $t>0$.
By the method in Section 4 in [Wa], similarly to Theorem 2.12 in [Wa], we get\\

\noindent{\bf Lemma 5.5} {\it When $2k\leq n$ and $X$ is
small, then for $f^j\in C^{\infty}_G(M)$,}
$${\rm lim}_{t\rightarrow 0}\phi_t\psi_t{{\bf {\rm ch}}^{2k}}(\sqrt{t}B,X)(f^0,\cdots ,
f^{2k})_J$$
$$=\frac{1}{(2k)!}{(2\pi\sqrt{-1})}^{-n/2}\int_{Z}f^0\wedge
df^1\wedge\cdots\wedge
df^{2k}\widehat{A}(F^Z_{\mathfrak{g}}(X))_J.\eqno(5.18)$$\\

Extending Theorem 3.9 to the family case, we have by Lemma 5.5 by\\

\noindent{\bf Theorem 5.6} {\it When $2k\leq n$ and $X$ is
small, then for $f^j\in C^{\infty}_G(M)$,}
$${\rm lim}_{t\rightarrow 0}\phi_t\psi_t{{\bf {\rm ch}}^{2k}}(\sqrt{t}B,X)(f^0,\cdots ,
f^{2k})$$
$$=\frac{1}{(2k)!}{(2\pi\sqrt{-1})}^{-n/2}\int_{Z}f^0\wedge
df^1\wedge\cdots\wedge
df^{2k}\widehat{A}(F^Z_{\mathfrak{g}}(X)).\eqno(5.19)$$\\

\indent By Theorems 5.4 and 5.6 and the following homotopy property, similarly to Corollary 3.11, we have\\

\noindent{\bf Corollary 5.7} {\it When $X$ is
small, then}
$${\rm Ch}_{e^{-X}}[{\rm Ind}(D_{{\rm Im}p,+,z})]
={(2\pi\sqrt{-1})}^{-n/2}\int_Z\widehat{A}(F^Z_{\mathfrak{g}}(X)){\rm Ch}({\rm Im}p).\eqno(5.20)$$\\

Let $B_t=\sqrt{t}\phi_t\psi_t(B^{-X})$ and ${\mathcal{F}}_t=B^2_t$. Then we have the homotopy formula:\\

\noindent{\bf Proposition 5.8}~{\it There is a cochain ${\rm ch}(B_t,\frac{dB_t}{dt},X)$ such that the following formula holds}
$$\frac{d{\rm ch}(B_t,X)}{dt}=-(b+B+d_{B_0}){\rm ch}(B_t,\frac{dB_t}{dt},X).\eqno(5.21)$$\\

\noindent{\bf Proof.} We know that $B_t$ is a superconnection on the infinite dimensional bundle $C^{\infty}(M,E)\rightarrow B_0$ which
we write ${\mathcal{E}}\rightarrow B_0$. Let $\widetilde{B_0}=B_0\times {\mathbb{R}}_+$, and $\widetilde{\mathcal{E}}$ be the superbundle
$\pi^*{\mathcal{E}}$ over $\widetilde{B_0}$, which is the pull-back to $\widetilde{B_0}$ of ${\mathcal{E}}$. Define a superconnection
$\widehat{B}$ on $\widetilde{\mathcal{E}}$ by the formula
$$(\widehat{B}\beta)(x,t)=(B_t\beta(\cdot,t))(x)+dt\wedge\frac{\partial\beta(x,t)}{\partial t}.\eqno(5.22)$$
The curvature $\widehat{\mathcal{F}}$ of $\widehat{B}$ is
$$\widehat{\mathcal{F}}={\mathcal{F}}_t-\frac{dB_t}{dt}\wedge dt,\eqno(5.23)$$
where ${\mathcal{F}}_t=B^2_t$ is the curvature of $B_t$. By the Duhamel principle, then
$$e^{-\widehat{\mathcal{F}}}=e^{-{\mathcal{F}}_t}+\left(\int^1_0e^{-u{\mathcal{F}}_t}\frac{dB_t}{dt}e^{-(1-u){\mathcal{F}}_t}du\right)\wedge dt.\eqno(5.24)$$
Let $f^0,\cdots,f^{2k}$ be in $C^{\infty}_G(M)$, then $[\widehat{B},f^j]=[B_t,f^j]$. We replace $K_t$ in (4.5) and (4.6) by the above $B_t$, then we define the cochain ${\rm ch}(B_t,\frac{dB_t}{dt},X)$.
 So by (5.24), we get on $C^{\infty}_G(M)$ that
$${\rm ch}(\widehat{B},X)={\rm  ch}(B_t,X)+{\rm ch}(B_t,\frac{dB_t}{dt},X)dt.\eqno(5.25)$$
Similarly to (5.15), we have
$$(b+B+d_{\widetilde{B_0}}){\rm ch}(\widehat{B},X)=0;~~(b+B+d_{B_0}){\rm ch}({B_t},X)=0.\eqno(5.26)$$
By (5.25) and (5.26), we get Proposition 5.8.~~$\Box$\\

 \noindent {\bf Acknowledgements.} This work
was supported by NSFC No.11271062 and NCET-13-0721. The author would like to thank Profs.
Weiping Zhang and Huitao Feng for very helpful suggestions and discussions. The author would
like to thank the referee for careful reading and helpful comments.\\

\noindent{\large \bf References}\\

\noindent[ABP]M. Atiyah, R. Bott, V. Patodi, On the heat equation and the index theorem, Invent. Math., 19 (1973), 279-330.

\noindent[Az]F. Azmi, The equivariant Dirac cyclic cocycle, Rocky Mountain J. Math. 30(2000), 1171-1206.

\noindent[BC]M. Benameur, A. Carey, Higher spectral flow and an entire bivariant JLO cocycle, J. K-theory, 11 (2013), 183-232.

\noindent[BGV]N. Berline, E. Getzler and M. Vergne, {\it Heat kernels
and Dirac operators,} Springer-Verlag, Berlin, 1992.

 \noindent[BV]N. Berline and M. Vergne, The
equivariant index and Kirillov character formula, Amer. J. Math.
107 (1985), 1159-1190.

\noindent[BV1]N. Berline and M. Vergne, A computation of the
equivariant index of the Dirac operators, Bull. Soc. Math. France
113 (1985), 305-345.

\noindent [Bi]J. M. Bismut, The infinitesimal Lefschetz formulas: a
heat equation proof, J. Func. Anal. 62 (1985), 435-457.

\noindent [BiF] J. M. Bismut and D. S. Freed, The analysis of
elliptic families II, Commun. Math. Phys. 107 (1986), 103-163.

\noindent[BF]J. Block and J. Fox, Asymptotic pseudodifferential operators and index theory, Contemp. Math., 105 (1990), 1-32.

\noindent [CH]S. Chern and X. Hu, Equivariant Chern character
for the invariant Dirac operators, Michigan Math. J. 44 (1997),
451-473.

\noindent[Co]A. Connes, Entire cyclic cohomology of Banach algebras and characters of
$\theta$-summable Fredholm module, K-Theory 1 (1988), 519-548.

\noindent[CM]A. Connes and H. Moscovici, Cyclic cohomology, the Novikov conjecture and hyperbolic groups,
Topology 29 (1990), 345-388.

\noindent [Do]H. Donnelly, Eta invariants for G-space, Indiana
Univ. Math. J. 27 (1978), 889-918.

\noindent [Fe]H. Feng, A note on the noncommutative Chern
character (in Chinese), Acta Math. Sinica 46 (2003), 57-64.

\noindent [FGV]H. Figueroa, J. Gracia-Bond\'{i}a and J.
V\'{a}rilly, {\it Elements of noncommutative geometry},
Birkh\"{a}user Boston, 2001.

\noindent[Ge1]E. Getzler, The odd Chern character in cyclic homology
and spectral flow, Topology 32 (1993), 489-507.

\noindent[Ge2]E. Getzler, Pseudodifferential operators on supermanifolds and the Atiyah-Singer index theorem,
Comm. Math. Phys. 92 (1983), no. 2, 163-178.

\noindent [GS]E. Getzler and A. Szenes, On the Chern character of
theta-summable Fredholm modules, J. Func. Anal. 84 (1989), 343-357.

\noindent[Go]S. Goette, Equivariant eta invariants and eta forms, J.
reine angew Math. 526 (2000), 181-236.

\noindent[JLO]A. Jaffe, A. Lesniewski and K. Osterwalder, Quantum K-theory: The Chern
character, Comm. Math. Phys. 118 (1988), 1-14.

\noindent[KL]S. Klimek and A. Lesniewski, Chern character in
equivariant entire cyclic cohomology, K-Theory 4 (1991), 219-226.

\noindent [LYZ]J. D. Lafferty, Y. L. Yu and W. P. Zhang, A
direct geometric proof of Lefschetz fixed point formulas, Trans.
AMS. 329 (1992), 571-583.

\noindent[Pf]M. Pflaum, The normal symbol on Riemannian manifolds. New York J. Math. 4 (1998), 97-125.

\noindent[PW]R. Ponge and H. Wang, Noncommutative Geometry and Conformal Geometry. II. Connes-Chern character and the local equivariant index theorem, arXiv:1411.3703.

\noindent[Si]B. Simon, Trace ideals and their applications, London Math. Soc., Lecture Note
35, 1979.

\noindent[Wa]Y. Wang, The noncommutative infinitesimal equivariant index formula, J. K-Theory, 14 (2014), 73-102.

\noindent[Wa1]Y. Wang, The equivariant noncommutative Atiyah-Patodi-Singer index theorem, K-Theory 37 (2006),
213-233.

\noindent [Wu] F. Wu, The Chern-Connes character for the Dirac
operators on manifolds with boundary, K-Theory 7 (1993), 145-174.

\noindent [Zh]W. P. Zhang, A note on equivariant eta
invariants, Proc. AMS. 108 (1990), 1121-1129.\\

 \indent{\it School of Mathematics and Statistics, Northeast Normal University, Changchun Jilin, 130024, China }\\
 \indent E-mail: {\it wangy581@nenu.edu.cn}\\

\end{document}